\documentclass{elsart}
\usepackage{amssymb,amsfonts}
\usepackage{graphicx}
\newcommand{\trace}{\mathop{\rm Tr}\nolimits}

\newcommand{\twomat}[4]{\left(\begin{array}{cc}#1&#2\\#3&#4\end{array}\right)}

\newcommand{\C}{{\mathbb{C}}}
\newcommand{\R}{{\mathbb{R}}}
\newcommand{\I}{{\mathbb{I}}}

\newcommand{\id}{\I}
\newcommand{\be}{\begin{equation}}
\newcommand{\ee}{\end{equation}}
\newcommand{\bea}{\begin{eqnarray}}
\newcommand{\eea}{\end{eqnarray}}
\newcommand{\beas}{\begin{eqnarray*}}
\newcommand{\eeas}{\end{eqnarray*}}

\newtheorem{theorem}{Theorem}

\newtheorem{corollary}{Corollary}

\newtheorem{proposition}{Proposition}

\newcount\minute
\newcount\hour
\def\currenttime{%
    \minute\time
    \hour\minute
    \divide\hour60
    \the\hour:\multiply\hour60\advance\minute-\hour\the\minute}
\begin{document}
\begin{frontmatter}
\title{Schur multiplier norms for Loewner matrices}
\author{Koenraad M.R.\ Audenaert}
\address{
Department of Mathematics,
Royal Holloway, University of London, \\
Egham TW20 0EX, United Kingdom}
\ead{koenraad.audenaert@rhul.ac.uk}
\date{\today, \currenttime}
\begin{keyword}
Commutator \sep Schur multiplier norm \sep Loewner matrix \sep concave function \sep convex function
%\MSC 15A60 XXX
\end{keyword}
%------------------------------------------------------------------ ABSTRACT
\begin{abstract}
We study upper bounds on the Schur multiplier norm of Loewner matrices for concave and convex functions.
These bounds then immediately lead to 
upper bounds on the ratio of Schatten $q$-norms of commutators
$||\;[A,f(B)]\;||_q/||\;[A,B]\;||_q$. 
We also consider operator monotone functions, for which sharper bounds are obtained.
\end{abstract}

\end{frontmatter}
%%%%%%%%%%%%%%%%%%%%%%%%%%%%%%%%%%%%%%%%%%%%%%%%%%%%%%%%%%%%%%%%%%%%%%%%%%%%%%%%%%%%%%%%%%%%%
%%%%%%%%%%%%%%%%%%%%%%%%%%%%%%%%%%%%%%%%%%%%%%%%%%%%%%%%%%%%%%%%%%%%%%%%%%%%%%%%%%%%%%%%%%%%%
%%%%%%%%%%%%%%%%%%%%%%%%%%%%%%%%%%%%%%%%%%%%%%%%%%%%%%%%%%%%%%%%%%%%%%%%%%%%%%%%%%%%%%%%%%%%%
\section{Introduction}
The main impetus behind the work presented in this paper was to find good upper bounds on the ratio
\be
||\;[A,f(B)]\;||_q/||\;[A,B]\;||_q\label{eq:ratio}
\ee
in terms of the spectrum of $B$,
where $A$ is a general $n\times n$ matrix, $B$ is a Hermitian $n\times n$ matrix,
$[A,B]$ is the commutator $[A,B]=AB-BA$,
$f$ is a given function $f:\R\to \R$, and $||\cdot||_q$ is the Schatten $q$-norm.

It is not hard to see that this problem immediately reduces to the problem of finding good upper bounds
on the Schur multiplier $q$-norms of the Loewner matrix of $f$ in $B$.
This will be shown in detail in Section \ref{sec:schur}.
The bulk of this paper is devoted to obtaining such bounds.

We will restrict attention to two classes of functions $f$: first the functions that are operator monotone on an
interval containing the spectrum of $B$ (see Theorem \ref{th:opmono} in Section \ref{sec:mono}),
and then in more generality the concave and convex functions
(see  Theorem \ref{th:conc} in Section \ref{sec:conc}).

Davies \cite{davies88} has considered a similar question for the function $x\mapsto|x|$ but he was looking for a universal
bound independent of  $B$.
He found that
for all Schatten norms except the trace norm and operator norm,
for all bounded operators $A$ and all self-adjoint operators $B$ in the Schatten $q$-class
$$
||\;[A,|B|]\;||_q/||\;[A,B]\;||_q \le 2(1+\gamma_q)
$$
where
$$
\gamma_q = c\min(q,q/(q-1))
$$
and $c\ge 1$ is an absolute constant.
For the trace norm and operator norm no finite constants $\gamma_1$, $\gamma_\infty$ exist.
For the Frobenius norm he found the sharper bound
$$
||\;[A,|B|]\;||_2 \le ||\;[A,B]\;||_2.
$$
In Section \ref{sec:appl} we apply the main theorem of Section \ref{sec:conc} to 
obtain bounds on $||\;[A,|B|]\;||/||\;[A,B]\;||$ in terms of the number of positive and negative eigenvalues of $B$.

%%%%%%%%%%%%%%%%%%%%%%%%%%%%%%%%%%%%%%%%%%%%%%%%%%%%%%%
\section{Schur multiplier norms\label{sec:schur}}
I begin by showing that finding a sharp upper bound on the ratio (\ref{eq:ratio})
amounts to computing the Schur multiplier norm (induced by a Schatten norm) of a Loewner matrix.

Let $L$ and $A$ be two matrices of the same dimension, then their entrywise product is denoted by $L\circ A$, i.e.\
$(L\circ A)_{ij} = L_{ij} A_{ij}$. This product is known alternatively as the Schur product (or Hadamard product).
The linear operator $S_L: A\mapsto L\circ A$ is called the Schur multiplier operator.
Any norm $|||\cdot|||$ on $A$ induces a norm on $S_L$, which we'll also denote by $|||S_L|||$, defined by
$$
|||S_L||| = \sup_A \frac{|||L\circ A|||}{|||A|||}.
$$
We will be interested in particular in the Schatten $q$-norms, which are defined as
$||A||_q = (\trace (A^* A)^{q/2})^{1/q}$. These include the trace norm ($q=1$), the Frobenius norm ($q=2$) and the operator norm
$||A||$ (the limit of $q\to\infty$).
The corresponding induced norms for $S_L$ are defined as
$$
||S_L||_q = \sup_A \frac{||L\circ A||_q}{||A||_q}.
$$

A basic property of any Schur multiplier norm is its self-duality.
If $|||\cdot|||_D$ is the dual norm of $|||\cdot|||$, then $|||S_L||| = |||S_L|||_D$.
In particular, $||S_L||_q = ||S_L||_{q'}$, where $1/q'=1-1/q$.
This can be proven easily using a standard duality argument. For example,
\cite{mathias93} gives a proof for the operator norm and its dual, the trace norm, but the proof works
for any other norm.
%Thus:
%$$
%||S_L|| := \max_{A} \frac{||S_L(A)||_\infty}{||A||_\infty} =
%||S_L||_1 := \max_{A} \frac{||S_L(A)||_1}{||A||_1}.
%$$

The importance of Schur multiplier norms for the problem considered in this paper
follows from the following proposition:
\begin{proposition}
Let $A$ be any matrix, and let $B$ be Hermitian with eigenvalues $b_i$.
Let $L$ be the Loewner matrix of $f$ at $B$:
$$
L_{ij} := \left\{
\begin{array}{ll}
\frac{f(b_i)-f(b_j)}{b_i-b_j},& b_i\neq b_j \\[2mm]
f'(b_i),& b_i=b_j.
\end{array}
\right.
$$
Then
\be
|||\;[A,f(B)]\;||| \le |||S_L|||\;\;|||\;[A,B]\;|||.
\ee
\end{proposition}
\textit{Proof.}
Working in the eigenbasis of $B$, the commutators can be expressed in terms of the Schur product
as follows:
$$
[A,B] = A\circ (b_i-b_j)_{i,j=1}^n,
\quad
[A,f(B)] = A\circ (f(b_i)-f(b_j))_{i,j=1}^n.
$$
%Clearly, the diagonals of both commutators are identically zero.
Consider now the Loewner matrix $L$ of the proposition.
It is easy to see that this can be expressed in terms of $L$ as
$$
[A,f(B)] = [A,B] \circ L = S_L([A,B]).
$$
Hence,
the norms of both commutators are related by
$$
|||\;[A,f(B)]\;||| \le |||S_L|||\;\;|||\;[A,B]\;|||.
$$
\qed

For the Schatten 2-norm (Frobenius norm), the induced Schur multiplier norm is easily calculated:
\bea
||S_L||_2 &=& \max_A \frac{||L\circ A||_2}{||A||_2} \nonumber\\
&=& \max_A \left(\frac{\sum_{i,j}|L_{ij}|^2 |A_{ij}|^2}{\sum_{i,j} |A_{ij}|^2}\right)^{1/2} \nonumber\\
&=& \max_{i,j} |L_{ij}|.\label{eq:S2}
\eea
Computing Schur multiplier norms for other norms than the 2-norm is in general very difficult,
and the fact that all entries of $L$ are in a certain range by no means implies that
$||S_L||$ should be in that range.
Indeed, when $L$ is upper triangular with all entries above the diagonal equal to 1, and all others 0, its
Schur multiplier norm is $O(\log n)$ \cite{cowen92}.
%For the considered classes of functions a sharp bound can be found that is dimension-independent, but this will on occasion
%require intricate arguments.

Using complex interpolation, bounds for general Schatten $q$-norms can be derived from bounds for the 1-norm and the 2-norm.
Indeed, by a direct application of Theorem 5.2 in Chapter 3 of \cite{gohberg70}, for any $1\le q\le2$ we have
\be
||S_L||_q \le ||S_L||_1^{2-q}  ||S_L||_2^{q-1}.\label{eq:interp}
\ee

%%%%%%%%%%%%%%%%%%%%%%%%%%%%%%%%%%%%%%%%%%%%%
\section{Operator monotone functions\label{sec:mono}}
The first and easiest class of functions treated here are the functions that are
operator monotone on a given interval $I$.
\begin{theorem}\label{th:opmono}
Let $f$ be an operator monotone function on the interval $I$.
Let $B$ be an $n\times n$ Hermitian matrix with spectrum in $I$.
Let $L$ be the Loewner matrix of $f$ at $B$.
Then, for all Schatten $q$-norms,
\be
||S_L||_q \le f'(\lambda_{\min}(B)).\label{eq:mono}
\ee
\end{theorem}
Note that here $f'$ is always non-negative over $I$.

\textit{Proof.}
If $f$ is operator monotone, then its Loewner matrix $L$ is a positive semidefinite matrix. By a theorem of Schur
(see \cite{bhatia07}, section 1.4),
$S_L$ is then a completely positive map and $||S_L||$ (and hence $||S_L||_1$) is equal to $\max_i L_{ii}$.
In the present case, this number is equal to $\max_i f'(b_i)$.
By the concavity of operator monotone functions, this maximum is equal to $f'(\min_i b_i)$.

For the Schatten 2-norm, we already found that
$||S_L||_2 = \max_{i,j} |L_{ij}|$.
Again, in the present case $\max_{i,j} |L_{ij}|=f'(\min_i b_i)$, which proves the inequality for the Frobenius norm.

Finally, using the complex interpolation bound (\ref{eq:interp}),
these two results imply that (\ref{eq:mono}) holds for all Schatten norms.
Indeed, for any $1\le q\le2$,
$$
||S_L||_{q'} = ||S_L||_q \le ||S_L||_1^{2-q}  ||S_L||_2^{q-1} \le f'(\min_i b_i).
$$
\qed
%%%%%%%%%%%%%%%%%%%%%%%%%%%%%%%%%%%%%%%%%%%%%%%%%%%%%%%%%%%%%%%%
\section{The numerical radius and its dual norm\label{sec:num}}
In this section I obtain an intermediary result needed in the next section, which may be of independent interest.

The numerical radius is defined as
$$
w(A) = \sup_x \frac{|\langle Ax|x\rangle|}{||x||^2}.
$$
This is a norm, and its dual norm is \cite{ando91}
$$
||Y||_{w^*} = \sup_X \frac{|\trace Y^*X|}{w(X)} = \sup_X \{|\trace Y^*X|:w(X)\le 1\},
$$
which I will call the $w^*$ norm here.
The unit ball of the $w^*$ norm
is the \emph{absolute convex hull} of the matrices of the form $xx^*$ with $x\in\C^n$ and $||x||=1$;
i.e.\ it is the set of matrices
$\sum_i \lambda_i x_i x_i^*$ for which $\sum_i |\lambda_i|\le 1$ and $||x_i||=1$.
This includes but is not limited to the normal matrices with trace norm not exceeding 1.

In general, the numerical radius never exceeds the spectral norm, $w(X)\le ||X||$.
Likewise, the $w^*$ norm is bounded below by the trace norm.
Indeed,
$$
||Y||_{w^*} = \sup_X \frac{|\trace Y^*X|}{w(X)}
\ge \sup_X \frac{|\trace Y^*X|}{||X||} = ||Y||_1.
$$

For normal matrices $X$, the numerical radius is equal to the spectral norm: $w(Y)=||Y||$.
Here we show that for normal matrices the $w^*$ norm is equal to the trace norm.
\begin{theorem}\label{th:wstar}
If $Y$ is normal, then
$||Y||_{w^*} = ||Y||_1$.
\end{theorem}
\textit{Proof.}
By a theorem of Ando \cite{ando73}, a matrix $X$ has numerical radius at most one if and only if
there exist contractions $W$ and $Z$, where $Z$ is Hermitian, such that
$$
X = (\id+Z)^{1/2} W (\id-Z)^{1/2}.
$$
The definition of the $w^*$ norm can therefore be rewritten as
\beas
||Y||_{w^*} &=& \sup_X \{|\trace Y^*X|:w(X)\le 1\} \\
&=& \sup_{W,Z} \{|\trace(Y^* (\id+Z)^{1/2} W (\id-Z)^{1/2})|: Z=Z^*, ||Z||\le 1, ||W||\le 1\} \\
&=& \sup_Z \left\{ \sup_W \left\{|\trace W(\id-Z)^{1/2} Y^* (\id+Z)^{1/2}|: ||W||\le 1\right\}: Z=Z^*, ||Z||\le 1\right\} \\
&=& \sup_Z \left\{ ||(\id-Z)^{1/2} Y^* (\id+Z)^{1/2}||_1: Z=Z^*, ||Z||\le 1\right\}.
\eeas
Since $Y$ is normal, it has a unitary spectral decomposition $Y=\sum_{j=1}^n \lambda_j u_j u_j^*$,
with $\{u_j\}_{j=1}^n$ an orthonormal basis of $\C^n$.
Hence,
$$
||(\id-Z)^{1/2} Y^* (\id+Z)^{1/2}||_1 \le
\sum_j |\lambda_j|\;||(\id-Z)^{1/2} u_j u_j^* (\id+Z)^{1/2}||_1.
$$
Noting that for any Hermitian contraction $Z$
\beas
||(\id-Z)^{1/2} u_j u_j^* (\id+Z)^{1/2}||_1 &=& \langle (\id-Z^2)^{1/2} u_j|u_j\rangle \\
&\le& ||(\id-Z^2)^{1/2}|| \le 1,
\eeas
we find
$$
||(\id-Z)^{1/2} Y^* (\id+Z)^{1/2}||_1 \le
\sum_j |\lambda_j| = ||Y||_1,
$$
and therefore
$$
||Y||_{w^*} \le ||Y||_1.
$$
\qed

\begin{corollary}\label{cor:SL}
For $n\times n$ Hermitian $L$,
$$
||S_L|| = \max_{x\in\C^n} \{||L\circ xx^*||_1:||x||\le 1\}.
$$
\end{corollary}
\textit{Proof.}
By Corollary 3 in \cite{ando91}, if $L$ is Hermitian, $||S_L||$ is equal to $||S_L||_w$,
the Schur multiplier norm of $S_L$ induced by the numerical radius:
$$
||S_L||_w := \sup_X \frac{w(L\circ X)}{w(X)}.
$$
By Lemma 1 in \cite{ando91}, $||S_L||_w\le 1$ if and only if for all vectors $x\in\C^n$,
$$
||\overline{L}\circ xx^*||_{w^*}\le ||x||^2.
$$
If $L$ is Hermitian, then so is $\overline{L}\circ xx^*$, so that by Theorem \ref{th:wstar},
$$
||\overline{L}\circ xx^*||_{w^*} = ||\overline{L}\circ xx^*||_1 = ||L\circ xx^*||_1.
$$
\qed
%%%%%%%%%%%%%%%%%%%%%%%%%%%%%%%%%%%%%%%%%%%%%%%%%%%%%%%%%%%%%%
\section{Concave and convex functions\label{sec:conc}}
In this section I consider the extension of Theorem \ref{th:opmono} to the concave and convex functions.
For these functions the Loewner matrix $L$ is no longer positive semidefinite in general.
However, it satisfies a number of monotonicity properties that will be useful in deriving upper bounds.
Let, as before, $B$ a Hermitian $n\times n$ matrix with eigenvalues $(b_j)_{j=1}^n$
sorted in non-decreasing order, $b_1\le b_2\le\ldots\le b_n$,
and denote the Loewner matrix of $f$ at $B$ by $L$.
The entries of $L$ are
$$
L_{ij} := \left\{
\begin{array}{ll}
\frac{f(b_i)-f(b_j)}{b_i-b_j},& b_i\neq b_j \\[2mm]
f'(b_i),& b_i=b_j.
\end{array}
\right.
$$
For concave $f$ and non-decreasing $b$, these elements satisfy the following relations:
$$
(R): \quad \left\{
\begin{array}{l}
L_{ij} = L_{ji};\\
\mbox{for } i\le j<k,\; L_{ij}\ge L_{ik};\\
\mbox{for } j<k\le i,\; L_{ji}\ge L_{ki}.
\end{array}
\right.
$$
As a consequence, for $i<j$, $L_{ii}\ge L_{jj}$, and for all $i$ and $j$,
$L_{11}\ge L_{ij}\ge L_{nn}$.

The case of the Frobenius norm is again very simple.
\begin{theorem}\label{prop:conc2}
Let $B$ be a Hermitian $n\times n$ matrix with eigenvalues $(b_j)_{j=1}^n$
sorted in non-decreasing order, $b_1\le b_2\le\ldots\le b_n$.
Let $f$ be a function that is concave or convex on the interval $[b_1,b_n]$.
Let $L$ be the Loewner matrix of $f$ at $B$.
Then
\be
||S_L||_2 \le \max(|f'(b_1)|, |f'(b_n)|).\label{eq:f2}
\ee
\end{theorem}
\textit{Proof.}
By (\ref{eq:S2}), the upper bound is given by $\max_{i,j} |L_{ij}|$.
For concave $f$, the properties (R) of $L$ imply that $\max_{i,j} |L_{ij}| = \max(|L_{11}|, |L_{nn}|)$.
Since $L_{ii}=f'(b_i)$ this proves inequality (\ref{eq:f2}).
For convex $f$, simply replace $f$ by $-f$ and note that both sides of the inequality are invariant under this sign change.
\qed

For the Schur multiplier trace norm (operator norm)
I start with a technical proposition about certain standardised monotonously increasing concave functions,
as the general case follows easily from this case.
\begin{proposition}\label{prop:conc}
Let $B$ be a Hermitian $n\times n$ matrix with eigenvalues $(b_j)_{j=1}^n$
sorted in non-decreasing order, $b_1\le b_2\le\ldots\le b_n$.
Let $g$ be a function that is concave on the interval $[b_1,b_n]$, and for which $g'(b_1)=1$ and $g'(b_n)=0$.
Let $K$ be the Loewner matrix of $g$ at $B$.
Then
\be
||S_K||_1 = ||S_K|| \le 1+\phi^{-1}\sum_{j=1}^{n}(1-g'(b_j)),
\ee
where $\phi$ is the Golden Ratio, $\phi=(1+\sqrt5)/2\approx 1.618$.
\end{proposition}
Note that the interpolation relation (\ref{eq:interp}) can again be used to obtain bounds for general Schatten norms.

\textit{Proof.}
The matrix $K$ satisfies conditions (R), and $K_{11}=1$ and $K_{nn}=0$.
From this I will derive an upper bound on $||S_K||$ in terms of the diagonal elements $k_j=K_{jj}$.

By Corollary \ref{cor:SL}, the Schur multiplier norm of $K$ can be characterised as
$$
||S_K||_1 = ||S_K|| = \max_{x\in \C^n} \{ ||K\circ(xx^*)||_1: ||x||=1\}.
$$
We can find an upper bound on the trace norm of any matrix $A$ by partitioning $A$
as the block matrix
$$
A = \twomat{B}{b}{b^T}{a},
$$
where $B$ is the upper left $(n-1)\times (n-1)$ submatrix of $A$, $a=A_{nn}$ and
$b$ is the $(n-1)$-dimensional
vector consisting of the first $(n-1)$ entries of the last column of $A$.
By a result of Bhatia and Kittaneh \cite{BK90}, the trace norm of $A$ can be bounded above by the sum
of the trace norms of the four blocks, i.e.\
$$
||A||_1 = ||B||_1 + 2||b||+|a|.
$$
When we apply this to the matrix $K\circ(xx^*)$,
we have $a=K_{nn}|x_n|^2=0$, $b_i=\overline x_nx_iK_{in}$
and $B_{ij}=K_{ij} x_i\overline x_j$, for $i,j=1,\ldots,n-1$.

Since the vector $x$ is normalised, the norm of the subvector of its first $n-1$ entries is
equal to $\sqrt{1-|x_n|^2}$.
Introducing the $(n-1)$-dimensional normalised vector $y$ with $y_i=x_i/\sqrt{1-|x_n|^2}$, for $i=1,\ldots,n-1$,
and partitioning $K$ conformally with $A$ as
$$
K = \twomat{Z}{u}{u^T}{0},
$$
we get
$b=\overline x_n \sqrt{1-|x_n|^2}\; (y\circ u)$
and $B=(1-|x_n|^2)\;(Z\circ(yy^*))$.
Hence
$$
||K\circ(xx^*)||_1 \le (1-|x_n|^2)\;||Z\circ(yy^*)||_1+2|x_n| \sqrt{1-|x_n|^2} \;||y\circ u||.
$$
As the maximisation over $x$ reduces to a maximisation over $|x_n|$ and over $y$, we obtain
$$
||S_K|| \le \max_{0\le x\le 1} (1-x^2)||S_Z||+2x \sqrt{1-x^2}\; \max_y \{||y\circ u||:||y||\le 1\}.
$$
The maximisation $\max_y \{||y\circ u||:||y||\le 1\}$ yields $\max_i u_i$,
which because of (R) is equal to $K_{1n}$ and therefore bounded above by 1.
Furthermore, substituting $a=||S_Z||$ and $x=\cos\theta$, the remaining maximisation is
$$
\max_{0\le \theta\le \pi/2} a(1-\cos2\theta)/2+\sin2\theta,
$$
which is the monotonously increasing function
$$
v(a):=a/2+\sqrt{1+(a/2)^2}.
$$
%Simple upper bounds are $v(a)\le 1+a$ and $v(a)^2\le 2+a^2$.
This gives our second relation:
\be
||S_K||\le v(||S_Z||).\label{eq:rec2}
\ee

Let us write $Z$ in terms of a matrix $K'$ with
upper left element 1 and lower right element 0: $Z=k_{n-1}J+(1-k_{n-1})K'$,
where $J$ is the $n\times n$ matrix with $J_{ij}=1$.
Note that $K'$ is a matrix that still obeys (R) but for which $k'_{n-1}=0$ and $k'_1=1$,
i.e.\ it has the same characteristics as the matrix $K$ we started out with.
The diagonal elements of $K'$ in terms of those of $K$ are given by
\be
k'_j:=\frac{k_j-k_{n-1}}{1-k_{n-1}}.\label{eq:rec3}
\ee
By convexity of the Schur multiplier norm and the fact that $||S_J||=1$, we have
$$
||S_Z||\le k_{n-1}+(1-k_{n-1})||S_{K'}||,
$$
so that, by (\ref{eq:rec2}),
\be
||S_K||\le v(k_{n-1}+(1-k_{n-1})||S_{K'}||)\label{eq:rec4}
\ee

The two relations (\ref{eq:rec3}) and (\ref{eq:rec4}) allow to find an
easily computable upper bound on $S_K$ via a recursion process.
This process stops after $n$ steps, as for a scalar $||S_a||=|a|$.
In the recursion, we need in succession the elements $k_{n-1},k'_{n-2}, k''_{n-3},\ldots,k^{(m)}_{n-m-1}$,
which I'll abbreviate by $a_m$, for $m=0,\ldots,n-2$.
Calculating it through, an explicit formula for the elements is
$$
a_0=k_{n-1}
$$
and, for $m=1,\ldots,n-2$,
$$
a_m=k^{(m)}_{n-m-1} = \frac{k_{n-m-1}-k_{n-m}}{1-k_{n-m}}.
$$
The last element in this sequence is (since $k_1=1$)
$$
a_{n-2} = \frac{k_{1}-k_{2}}{1-k_{2}} = 1.
$$
Then, denoting $||S_{K^{(m)}}||$ by $s_m$,
$$
s_m \le v(a_m +(1-a_m)s_{m+1}),\quad
s_{n-2}=1.
$$
Defining $t_m=s_m-1$ and
$$
b_m = 1-a_m = \frac{1-k_{n-m-1}}{1-k_{n-m}},
$$
we have
$$
t_m \le v(1+b_m t_{m+1})-1,\quad
t_{n-2}=0.
$$
It is easily verified that $v(1+x)-1\le 1/\phi+x$, where $\phi$ is the Golden Ratio.
Thus
$$
t_m \le b_m t_{m+1}+1/\phi,\quad
t_{n-2}=0,
$$
whence
$$
t_0 \le \phi^{-1}(1 + b_0+b_0b_1+\ldots+b_0b_1\cdots b_{n-3}).
$$
It is immediately checked that
$b_0b_1\cdots b_j = 1-k_{n-j-1}$, for $j=0,\ldots,n-3$
and $k_1=1$, $k_n=0$,
so that
$$
t_0 \le \phi^{-1}\sum_{j=1}^n(1-k_j).
$$
This finally yields
$||S_K|| \le s_0 \le 1+ \phi^{-1}\sum_{j=1}^n(1-k_j)$.
As $K_{ii}=g'(b_i)$, the inequality of the proposition follows.
\qed

%%%%%%%%%%%%%%%%%%%%%%%%%%%%%%%%%%%%%%%%%%%%%%%%%%%%%%%%%%%%%%%
%%%%%%%%%%%%%%%%%%%%%%%%%%%%%%%%%%%%%%%%%%%%%%%%%%%%%%%%%%%%%%%
\begin{corollary}\label{cor:conc}
Let $B$ be a Hermitian $n\times n$ matrix with eigenvalues $(b_j)_{j=1}^n$
sorted in non-decreasing order, $b_1\le b_2\le\ldots\le b_n$.
Let $h$ be a function that is concave on the interval $[b_1,b_n]$, and for which $h'(b_1)=0$ and $h'(b_n)=-1$.
Let $K$ be the Loewner matrix of $h$ at $B$.
Then
\be
||S_K|| \le 1+\phi^{-1}\sum_{j=1}^{n}(1+h'(b_j)).
\ee
\end{corollary}
\textit{Proof.}
This follows immediately from Proposition \ref{prop:conc}
with the matrix $B$ replaced by $B'=b_1+b_n-B$ 
and defining $h(x)=g(b_1+b_n-x)$, so that $h'(b_j)=-g'(b_1+b_n-b_j)=-g'(b'_j)$.
\qed

We can now state and prove the main result of this paper.
\begin{theorem}\label{th:conc}
Let $B$ be a Hermitian $n\times n$ matrix with eigenvalues $(b_j)_{j=1}^n$
sorted in non-decreasing order, $b_1\le b_2\le\ldots\le b_n$.
Let $f$ be a function that is concave on the interval $[b_1,b_n]$.
Let $L$ be the Loewner matrix of $f$ at $B$.
Then
\beas
||S_L|| &\le& (\alpha-\beta)+\min\Big(|\beta|+\phi^{-1}\;\sum_{j=1}^{n}(\alpha-f'(b_j)),\\
&&                               |\alpha|+\phi^{-1}\;\sum_{j=1}^{n}(f'(b_j)-\beta)\Big),
\eeas
where $\alpha=f'(b_1)$ and $\beta=f'(b_n)$.
For any function that is convex on the interval $[b_1,b_n]$,
\beas
||S_L|| &\le& (\beta-\alpha)+\min\Big(|\beta|+\phi^{-1}\;\sum_{j=1}^{n}(f'(b_j)-\alpha),\\
                               |\alpha|+\phi^{-1}\;\sum_{j=1}^{n}(\beta-f'(b_j))\Big).
\eeas
\end{theorem}
\textit{Proof.}
General concave functions $f$ can be mapped to the standardised functions $g$ and $h$ of Proposition \ref{prop:conc}
and Corollary \ref{cor:conc}.
Note that
$$
\alpha:=f'(b_1) \ge f'(b_j)\ge f'(b_n)=:\beta.
$$

First we write
$$
f(x) = \beta x+(\alpha-\beta)g(x).
$$
Then
$$
(\alpha-\beta)g'(x) = f'(x)-\beta.
$$
Letting $L$ and $K$ be the Loewner matrices of $f$ and $g$, respectively, at $B$,
$$
L=\beta J+(\alpha-\beta)K,
$$
where $J$ is the matrix all of whose entries are 1.
As $||S_J||=1$,
\beas
||S_L||
&\le& |\beta| + (\alpha-\beta)||S_K|| \\
&\le& |\beta| + (\alpha-\beta)\left(1+\phi^{-1}\sum_{j=1}^{n}(1-g'(b_j))\right) \\
&=& |\beta| + (\alpha-\beta) + \phi^{-1}\sum_{j=1}^{n}((\alpha-\beta)-(f'(b_j)-\beta)) \\
&=& |\beta| + (\alpha-\beta) + \phi^{-1}\sum_{j=1}^{n}(\alpha-f'(b_j)).
\eeas

We can also write
$$
f(x) = \alpha x+(\alpha-\beta)h(x).
$$
and obtain in a similar way
$$
||S_L|| \le
 |\alpha| + (\alpha-\beta) + \phi^{-1}\sum_{j=1}^{n}(f'(b_j)-\beta).
$$

Taking the minimum of both bounds yields the bound of the corollary.

For convex $f$ we just replace $f$ by $-f$ and apply the result for concave functions.
Since now $\alpha:=f'(b_1) \le f'(b_j)\le f'(b_n)=:\beta$, an appropriate sign change has to be applied to the bound.
\qed

When the spectrum of $B$ is not known, but it is known that $b_1\le B\le b_n$,
weaker bounds follow readily from this Theorem:
\begin{corollary}
Let $B$ be a Hermitian $n\times n$ matrix bounded as  $b_1\le B\le b_n$.
Let $f$ be a function that is either concave or convex on the interval $[b_1,b_n]$.
Let $L$ be the Loewner matrix of $f$ at $B$.
Then
$$
||S_L|| \le |\alpha-\beta|(1+(n-1)\phi^{-1})+\min\left(|\beta|,|\alpha|\right),
$$
where $\alpha=f'(b_1)$ and $\beta=f'(b_n)$.
\end{corollary}
%%%%%%%%%%%%%%%%%%%%%%%%%%%%%%%%%%%%%%%%%%%%%%%%%%%%%%%%%%%%%%%%%%%%%%%%%%%%%%%%
\section{Examples\label{sec:appl}}
As a first application, we consider the function $f(x) = |x|$.
\begin{theorem}
Let $B$ be a Hermitian $n\times n$ matrix with $r$ positive eigenvalues.
Let $L$ be the Loewner matrix of the function $f(x)=|x|$ at $B$.
Then, for $1\le r<n$,
$$
||S_L|| \le 3+2\phi^{-1}\min(r,n-r).
$$
If $r$ is $0$ or $n$, $||S_L||$ is 1.
\end{theorem}
\textit{Proof.}
For $1\le r<n$,  $\alpha=f'(b_1)=1$, $\beta=f'(b_n)=-1$, $f'(b_j)=1$ for $r$ values of $j$, and
$f'(b_j)=-1$ for $n-r$ values of $j$.
The bound follows by simple calculation.
\qed

Since the bounds only depend on the diagonal elements of the Loewner matrix, they are not expected to be
sharp for specific functions. For the absolute value function, for example, it is known that in the $d=2$ case
the norm ratio lies between the values $1$ and $\sqrt{2}$, whereas the Theorem gives the bound $3+2/\phi$ for $r=1$.

\bigskip

For our second example, consider the following corollary of the main theorem.
Let $C$ be a Hermitian matrix with spectrum $c_1\le c_2\le \ldots \le c_n$.
By putting $B=g(C)$ and $h=f\circ g$, we find:
\begin{corollary}
For all $n\times n$ matrices $A$ and for any monotonously increasing function $g$ and any function $h$ such that
$f=h\circ g^{-1}$ is concave,
\beas
\frac{||\;[A,h(C)]\;||_1}{||\;[A,g(C)]\;||_1}
&\le& (\alpha-\beta)+\min\Big(|\beta|+\phi^{-1}\;\sum_{j=1}^{n}(\alpha-f'(g(c_j))),\\
&&                               |\alpha|+\phi^{-1}\;\sum_{j=1}^{n}(f'(g(c_j))-\beta)\Big),
\eeas
where $\alpha=f'(g(c_1))$ and $\beta=f'(g(c_n))$.
\end{corollary}

Consider the functions
$h(x)=\log x$ and $g(x)=\log(x)-\log(1-x)$.
Thus, $f(x)=x-\log(1+e^x)$, which is monotonously increasing and concave,
and $(f'\circ g)(x)=1-x$.
The bound of the corollary then simplifies to
$$
c_n-c_1+\min\left(1-c_n+\phi^{-1}(1-nc_1),1-c_1+\phi^{-1}(nc_n-1)\right),
$$
As $c_1\ge 0$, this quantity is bounded above by $1+\phi^{-1} = \phi$. We have therefore proven:
\begin{corollary}
For any $A$ and for any positive semidefinite $C$ with $\trace C=1$,
\be
||\;[A,\log(C)]\;||_1 \le \phi\;\;||\;[A,\log(C)-\log(\id-C)]\;||_1.
\ee
\end{corollary}

%------------------------------------------------------------- BIBLIOGRAPHY
%%%\section*{Acknowledgments}

%%%%%%%%%%%%%%%%%%%%%%%%%%%%%%%%%%%%%%%%%%%%%%%%%%%%%%%%%%%%%%%%%%%

\begin{thebibliography}{9}
%\bibitem{bhatia} R.~Bhatia, \textit{Matrix Analysis}, Springer, Heidelberg (1997).
\bibitem{ando73} T.~Ando, ``On the structure of operators with numerical radius one'',
Acta Sci.\ Math.\ (Szeged) \textbf{34}, 11--15 (1973).
\bibitem{ando91} T.~Ando and K.~Okubo, ``Induced norms of the Schur multiplier operator'',
Linear Algebra Appl.\ \textbf{147}, 181--199 (1991).
\bibitem{cowen92} J.R.~Angelos, C.C.~Cowen and S.K.~Narayan,
``Triangular truncation and finding the norm of a Hadamard multiplier'',
Linear Algebra Appl.\ \textbf{170}, 117--135 (1992).
\bibitem{bhatia07} R.~Bhatia, \textit{Positive Definite Matrices}, Princeton University Press, Princeton (2007).
\bibitem{BK90} R.~Bhatia and F.~Kittaneh,
``Norm inequalities for partitioned operators and an application'',
Math.\ Ann.\ \textbf{287}, 719--726 (1990).
%\bibitem{BK97} R.~Bhatia and F.~Kittaneh,
%``Some Inequalities for Norms of Commutators'',
%SIAM.\ J.\ Matrix Anal.\ \& Appl., \textbf{18}(1), 258--263 (1997).
%\bibitem{cowen96} C.C.~Cowen, P.A.~Ferguson, D.K.~Jackman, E.A.~Sexauer, C.~Vogt and H.J.~Woolf,
%``Finding norms of Hadamard multipliers'',
%Linear Algebra Appl.\ \textbf{247}, 217--235 (1996).
\bibitem{davies88} E.B.~Davies, ``Lipschitz continuity of functions of operators in the Schatten classes'',
J.\ London Math.\ Soc.\ (2) \textbf{37}(1), 148--157 (1988).
\bibitem{gohberg70} I.C.~Gohberg and M.G.~Krein, \textit{Theory and applications of Volterra operators in Hilbert space},
English translation, AMS, Providence (1970).
\bibitem{mathias93} R.~Mathias, ``The Hadamard Operator Norm Of A Circulant And Applications'',
SIAM J.\ Matrix Anal.\ Appl.\ \textbf{14}(4), 1152--1167 (1993).
%\bibitem{mathias93b} R.~Mathias, ``Matrix completions, norms, and Hadamard products'',
%Proc.\ AMS \textbf{117}(4), 905--918 (1993).
%\bibitem{olevskii94} V.~Olevskii and M.~Solomyak, ``An estimate for Schur multipliers in $S_p$-classes'',
%Linear Algebra Appl.\ \textbf{208/209}, 57--64 (1994).
%\bibitem{HM} R.A.~Horn and R.~Mathias, ``Cauchy-Schwarz inequalities associated with positive
%semidefinite matrices'', Linear Algebra Appl.\ \textbf{142}, 63--82 (1990).
\end{thebibliography}
\end{document}